\documentclass[a4paper,10pt]{article}
\usepackage[T1,T2A]{fontenc}
\usepackage[cp1251]{inputenc}
\usepackage[english]{babel}
\usepackage{amsmath}
\usepackage{amsfonts}
\usepackage{amssymb}
\usepackage{amsthm}
\usepackage{euscript}

\newtheorem{definition}{Definition}

\newtheorem{theorem}{Theorem}
\newtheorem{corollary}{Corollary}
\newtheorem{remark}{Remark}

\begin{document}
\author{P.~B.~Zatitskiy\thanks{Supported by Chebyshev Laboratory (SPbSU), RF Government grant 11.G34.31.0026, by JSC ``Gazprom Neft'', by the RFBR grants 13-01-12422 ofi\_m2 and 14-01-00373\_A, by President of Russia grant for young researchers MK-6133.2013.1. The author acknowledge Saint-Petersburg State University for a research grant 6.38.223.2014}}
\title{On scaling entropy sequence of dynamical system}

\maketitle

\begin{abstract}
We present a series of statements about scaling entropy, a metric invariant 
of dynamical systems proposed by A.~M.~Vershik in the late 90's (see~\cite{Versh10},~\cite{VershGorb07}, and~\cite{Versh11}). We show that it is a metric invariant indeed and calculate it in several special cases.
\end{abstract}

All measure spaces considered below are standard probability spaces (Lebesgue--Rokhlin spaces). Probability spaces with continuous measures are of major interest of this article, however all the definitions deal with general Lebesgue spaces (i.e. its measures may have 
atoms). We omit the sigma-algebra notation when we talk about probability spaces.

\section{Definition and properties of scaling entropy sequence}

Recall the definition of admissible semimetric and its $\varepsilon$-entropy (see, e.g.,~\cite{VZP},~\cite{ZP}):
\begin{definition}
Semimetric $\rho$ on a standard probability space $(X,\mu)$ is called measurable if it is measurable on $(X^2,\mu^2)$ as a function of two variables, and admissible if it is separable on some subset $X_1 \subset X$ of full measure. In this latter case we call the triple $(X,\mu,\rho)$ an admissible semimetric triple.
\end{definition}

The following definition goes back to A.~N.~Kolmogorov.

\begin{definition}
The function $\varepsilon \mapsto \mathbb{H}_{\varepsilon}(X,\mu,\rho)$ is called an epsilon entropy of a semimetric triple $(X, \mu, \rho)$ if it is defined for $\varepsilon>0$ as follows.   
Let $k\in \mathbb{N}$ be the least natural number, such that the space $X$ can be partitioned into disjoint measurable subsets $X_0,X_1, \dots, X_k$, such that $\mu(X_0)<\varepsilon$
and for any $j=1,\dots,k$ the diameter of the set $X_j$ in semimetric $\rho$ 
is less than $\varepsilon$.
Define the $\varepsilon$-entropy of the triple $(X,\mu, \rho)$
as $$\mathbb{H}_{\varepsilon}(X,\mu,\rho)=\log_2 k.$$
If one can not find such a number $k$, then we define $\mathbb{H}_{\varepsilon}(X,\mu,\rho)=+\infty.$
\end{definition}

Recall that a measurable semimetric $\rho$ on a standard probability space $(X,\mu)$ is admissible if and only if $\mathbb{H}_{\varepsilon}(X,\mu,\rho)<+\infty$ for any 
$\varepsilon>0$ (see~\cite{VZP} for details).

In what follows we assume $T$ to be an ergodic automorphism of a standard probability space $(X,\mu)$. Let $\rho$ be a measurable semimetric on $(X,\mu)$. Clearly, the
shifts $T^k\rho(x,y)= \rho(T^kx, T^ky)$ 
are also measurable semimetrics on $(X,\mu)$, moreover, for any $\varepsilon>0$ one has
$$
\mathbb{H}_\varepsilon(X,\mu,\rho)=\mathbb{H}_\varepsilon(X,\mu,T^k\rho).
$$
Define finite averages of a semimetric $\rho$ with respect to $T$ by the formula 
$
T^n_{av}\rho= \frac{1}{n}\sum\limits_{k=0}^{n-1}T^k\rho.
$

The following generating property is an important dynamical property of semimetrics:
\begin{definition}
We say that a measurable semimetric $\rho$ on $(X,\mu)$ is (two-sided) generating for a metric dynamical system $(X,\mu, T)$ if its shifts $T^k\rho, k\in \mathbb{Z},$ separate points mod zero, i.e. there exists some measurable subset $X'\subset X$ of full measure, such that for any two different $x,y \in X'$ one can find $k \in \mathbb{Z}$, such that $T^k\rho(x,y)>0$. 
\end{definition}
Note that any admissible metric is a generating semimetric for any automorphism $T$ because it separates points itself.

Recall the definition given by A.~M.~Vershik in~\cite{Versh10},~\cite{Versh11},~\cite{VershGorb07}.
\begin{definition}\label{scaldefsem}
Let $(X,\mu, T)$ be a metric dynamical system and $\rho$ be an admissible semimetric on $(X,\mu)$. A sequence $\{h_n\}$ of positive numbers is called a scaling sequence for the semimetric $\rho$, if for any $\varepsilon>0$ one has
\begin{equation}\label{upbound}
 \limsup\limits_n \frac{\mathbb{H}_\varepsilon(X,\mu, T^n_{av}\rho)}{h_n}<+\infty,
\end{equation}
and for $\varepsilon>0$ small enough one has
\begin{equation}\label{lowbound}
0< \liminf\limits_n \frac{\mathbb{H}_\varepsilon(X,\mu,T^n_{av}\rho)}{h_n}.
\end{equation}
We denote by $\EuScript{H}(X,\mu,T,\rho)$ the class of all scaling 
sequences for the semimetric $\rho$.
\end{definition} 

\begin{remark}
If $\rho$ is some admissible semimetric, $h=\{h_n\}\in \EuScript{H}(X,\mu,T,\rho)$, and $h'=\{h_n'\}$ is a sequence of positive numbers, then $h' \in \EuScript{H}(X,\mu,T,\rho)$ if and only if $0<\liminf \frac{h_n'}{h_n}\leq\limsup \frac{h_n'}{h_n}<\infty$.
\end{remark}

The following theorem confirms the hypothesis 
of Vershik~\cite{Versh11} about independence of a scaling 
sequence from a semimetric. It was proved in~\cite{VZP} (see Theorem~\ref{pointspect} below) for automorphisms with purely point spectrum. 
\begin{theorem}\label{scaling}
Let $\rho, \tilde{\rho} \in \mathcal{A}dm(X,\mu)$ be two
summable admissible generating semimetrics on~$(X,\mu)$. Then $\EuScript{H}(X,\mu,T,\rho)=\EuScript{H}(X,\mu,T,\tilde{\rho})$.
\end{theorem}

Theorem~\ref{scaling} leads to the following definition.
\begin{definition}\label{scaldef}
A sequence $h=\{h_n\}$ of positive numbers is called a scaling entropy sequence of a metric dynamical system $(X,\mu,T)$ if $h\in \EuScript{H}(X,\mu,T,\rho)$ for some (and then for any) summable admissible generating semimetric $\rho \in \mathcal{A}dm(X,\mu)$. We will denote by $\EuScript{H}(X,\mu,T)$ the class of all scaling entropy sequences for the dynamical system $\EuScript{H}(X,\mu,T)$.
\end{definition}

\begin{corollary}
The class of scaling entropy sequences is a metric invariant of ergodic dynamical systems.
\end{corollary}

It is an interesting question whether a scaling entropy sequence 
exists for any dynamical system (i.e. the set $\EuScript{H}(X,\mu,T)$ is nonempty).

The class of scaling entropy sequences is monotone in the following sense:
\begin{corollary}
Suppose that a metric dynamical system $(X_2,\mu_2,T_2)$ is a factor of a system $(X_1,\mu_1,T_1)$, and $h^1=\{h^1_n\}$, $h^2=\{h^2_n\}$ are corresponding scaling entropy sequences. Then $h^2_n=O(h^1_n)$, $n \to +\infty$.
\end{corollary}

\section{Examples}

\subsection{Relation to Kolmogorov entropy}
The following theorem clarifies the relation between 
scaling entropy sequences and Kolmogorov entropy.
\begin{theorem}\label{scalevskolm}
Let $T$ be an ergodic automorphism of a standard probability space $(X,\mu)$. Then:
\begin{itemize}
\item[1)] if Kolmogorov entropy $h_\mu(T)$ is finite and positive, then the sequence $h_n=n$ is a scaling entropy sequence of a system $(X,\mu,T)$;
\item[2)] if $h_n$ is a scaling entropy sequence of a system $(X,\mu,T)$, then $h_\mu(T)=0$ if and only if $h_n=o(n)$, $n \to +\infty$.
\end{itemize} 
\end{theorem}

It is not difficult to calculate scaling entropy sequences of Bernoulli shifts.
Let $(A,\nu)$ be a standard probability space, $X=A^\mathbb{Z}$ be the space of two-sided sequences, $\mu = \nu^\mathbb{Z}$ be the product measure on $X$, and $T\colon X \to X$ be the left shift. 
\begin{theorem}\label{scalebern}
If the measure $\nu$ does not have atoms, then the sequence $h_n=n$ is a scaling entropy sequence of the Bernoulli shift $(X,\mu,T)$. 
\end{theorem}
If the measure $\nu$ is continuous, then the entropy of the space $(A,\nu)$ is infinite. In this case according to the Kolmogorov theorem (see~\cite{Kolm58,Kolm59, Sinai59}) Kolmogorov entropy $h_\mu(T)$ is infinite. 
Therefore the converse of Theorem~\ref{scalevskolm} first statement does not hold
in general.

\subsection{Automorphisms with pure point spectrum}
The following characterization of automorphisms with pure point spectrum in terms of scaling entropy sequences was proved in~\cite{VZP}.
\begin{theorem}
\label{pointspect}
Let $T$ be an ergodic automorphism of a standard probability space $(X,\mu)$. Then the following are equivalent:
\begin{itemize}
\item[1)] the spectrum of the system $(X,\mu,T)$ is purely point;
\item[2)] the set $\EuScript{H}(X,\mu,T)$ consists of all positive bounded separated from zero sequences.
\end{itemize}
\end{theorem}

We compare this result with another criterion of the pure pointness of the dynamical system spectrum due to A.~G.~Kushnirenko (see~\cite{Kush67}). Recall the definition of $A$-entropy (sequential entropy) of an ergodic automorphism $T$ of a standard probability space $(X,\mu)$. Suppose that $A=\{a_n\}$ is an increasing sequence of natural numbers and $\xi$ is a measurable partition of $(X,\mu)$ with 
finite entropy $H(\xi)<+\infty$ ($H$ stands for entropy of measurable partition). Define $h_A(T,\xi)=\limsup\limits_{n}\frac{H(\cap_{i=1}^n T^{a_i}\xi)}{n}$ and $A$-entropy $h_A(T)=\sup\limits_\xi h_A(T,\xi)$ of automorphism $T$, where supremum is taken over all partitions with finite entropy. In~\cite{Kush67} Kushnirenko proved that an automorphism $T$ has purely point spectrum if and only if $h_A(T)=0$ for any sequence $A$. We compare this result with Theorem~\ref{pointspect} and obtain that if $h_A(T)>0$, then the scaling entropy sequence increases. The following statement gives a quantitative estimate of this remark.

\begin{theorem}\label{kush}
If $A=\{a_n\}$, $h_A(T)>0$, and $h=\{h_n\}$ is the scaling entropy sequence of an automorphism $T$, then $\liminf \frac{h_{a_n}}{\log\log n}>0$.
\end{theorem}

\subsection{Substitutions of constant length}
Now we calculate the scaling entropy sequences of substitutional dynamical systems of constant length. We refer the reader to the monograph~\cite{Queffelec} for detailed information concerning substitutional dynamical systems.

Let $\xi$ be a substitution of length $q$ on an alphabet $A$ (i.e. it is a map from $A$ to $A^q$). We assume that for some symbol $\alpha \in A$ the first symbol of the word $\xi(\alpha)$ is $\alpha$ again. In such a case one can find an infinite word $u \in A^{\mathbb{N}_0}$ starting with $\alpha$ such that $\xi(u)=u$. Let $T$ be the left shift on $A^{\mathbb{N}_{0}}$, and $X_\xi$ be the closure of the sequence $T^n u, n \in \mathbb{N},$ in the product topology. If the substitution $\xi$ is primitive, then the topological dynamical system $(X_\xi,T)$ is minimal and unique ergodic with unique invariant measure $\mu^\xi$. In such a manner we associate the metric dynamical system $(X_\xi, \mu^\xi,T)$ with the primitive substitution $\xi$. It turns out that the scaling entropy sequence of this dynamical system can be expressed in terms of two combinatorial parameters of the substitution. First of them, the height $h(\xi)$, can be defined for example as the greatest number $k$ relatively prime to $q$ such that if $u_n=u_0$ for some $n$, then $k \mid n$. Second, the column number $c(\xi)$, is defined as $\min\{|\{\xi^n(a)_k\colon a \in A\}|\colon n \in \mathbb{N}, 0\leq k < q^n \}$ 
(note that this definition differs a little from the classical one given in~\cite{Dekking77} and~\cite{Queffelec}).

\begin{theorem}\label{scalsubst}
Let $\xi$ be an injective primitive constant length substitution on an alphabet $A$, $|A|>1$. Then the sequence $h_n=1+(\mathrm{c}(\xi)-h(\xi))\log n$ is a scaling entropy sequence of the substitutional dynamical system $(X_\xi, \mu^\xi,T)$.
\end{theorem}

The idea of the proof of this theorem is to study $\xi$-invariant semimetrics on an alphabet $A$.

Theorems~\ref{scalsubst}~and~\ref{pointspect} imply the following corollary.
\begin{corollary}\label{substpoint}
Let $\xi$ be an injective primitive constant length substitution on an alphabet $A$, $|A|>1$. Then the spectrum of the corresponding dynamical system $(X_\xi, \mu^\xi,T)$ is 
purely point if and only if $\mathrm{c}(\xi)=h(\xi)$.
\end{corollary}

Note that the similar criterion was first obtained in~\cite{Kamae72} and~\cite{Dekking77}, but it was formulated a little bit different there. The difference is as follows. If the height $h(\xi)$ is bigger than one, then $\mathrm{c}(\xi)$ is defined via some modified substitution. The result of Corollary~\ref{substpoint} coincides with their criterion for the case $h(\xi)=1$: the spectrum is pure point if and only if $\mathrm{c}(\xi)=1$. But if $h(\xi)>1$, then criterion~\ref{substpoint} (unlike the Kamae-Dekking criterion) does not require to consider modified substitutions.

\bigskip
\bigskip

The author thanks F.~V.~Petrov and A.~M.~Vershik for a lot of discussions, advises and help.

Chebyshev Laboratory, St. Petersburg State University, 199178, Russia, Saint Petersburg, 14th Line, 29b;
PDMI RAS, 191023, Russia, St. Petersburg, Fontanka, 27.   
e-mail: paxa239@yandex.ru

\end{document}